\documentclass{article}
\usepackage[english,russian,ukrainian]{babel}
\usepackage{latexsym,amssymb}
\usepackage[cp1251]{inputenc}
\usepackage{amsfonts,amssymb}
\usepackage{euscript}
\newcounter{theorem} 
\newcounter{lemma} 

\begin{document}

\vspace*{7mm}

\Large







\noindent {\bf  Vitalii Shpakivskyi}\\

\noindent {\bf Construction of an infinite-dimensional family of exact solutions of a three-dimensional biharmonic equation by the hypercomplex method}

\vspace{7mm}

An infinite-dimensional family of exact solutions of a three-dimensional biharmonic equation was constructed by the hypercomplex method.

\Large

\section*{Preface}

In 2019, at the 12th ISAAC Congress (Aveiro, Portugal), I gave a report on the topic: Hypercomplex method for solving linear PDEs. 
After the report, Professor Yuri Grigoriev approached me and asked: can I write out the solutions of the three-dimensional biharmonic equation?
I answered: yes, I can. Then he asked me to write a separate paper on the solutions of the three-dimensional biharmonic equation and give him a copy.
I said: ok. But only now I am fulfilling my promise. It is a pity that Professor Yuri Grigoriev did not read it in his lifetime.
Only now did I understand that he had been studying the spatial theory of elasticity throughout his life, in particular, the three-dimensional biharmonic equation...

\section{Introduction}

Consider the three-dimensional biharmonic equation
 \begin{equation}\label{3-bih}
\Delta^2U(x,y,z):=\left(\frac{\partial^2 }{\partial
x^2}+\frac{\partial^2 } {\partial y^2}+\frac{\partial^2 } {\partial y^2}\right)^2U(x,y,z)=0
\end{equation}
or
$$
 \frac{\partial^4 U }{\partial
x^4}+\frac{\partial^4 U }{\partial
y^4}+\frac{\partial^4 U }{\partial
z^4}+2\frac{\partial^4 U} {\partial x^2\partial y^2}+2\frac{\partial^4 U} {\partial x^2 \partial z^2}+2\frac{\partial^4 U} {\partial y^2 \partial z^2}=0.
$$

Not many authors have tried to solve the three-dimensional biharmonic equation.
The reason, of necessity, is that three-dimensional problems require large computing power. 
 Such computing power has only recently begun to become available for research.

However, in paper \cite{Grigoriev-2014}, Yu. Grigoriev considers equation (\ref{3-bih}) in  star domains
relative to the origin of coordinates.
In such domains, he gives a representation of the general solution of  equation (\ref{3-bih}) in spherical  coordinates.
We will also note the paper \cite{Grigoriev-85} in which the theory of regular quaternion functions systematically were
 used for a solution of a three-dimensional problem of the elasticity theory.

In the paper  \cite{Lin-94}, were are concerned with energy decay estimates for
solutions of the biharmonic equation in a semi-infinite cylinder with nonzero boundary conditions on the finite end.
In the work  \cite{Altas-2002}, it is consider several finite-difference approximations for the three-dimensional biharmonic equation. 

In this paper we construct an infinite-dimensional family of exact solutions of equation (\ref{3-bih}) via holomorphic functions of a complex variable.

\section{The method}

This method originates from the works of P.\,W. Ketchum \cite{Ketchum-28,Ketchum-29-3-6,Ketchum-32-3-6}, 
and develop independently in the works of I.\,P. Mel'nichenko \cite{Mel'nichenko75,Melnichenko03}, 
S.\,A. Plaksa  \cite{Plaksa12} and ours.
 The final formation of the method was done in our papers \cite{Shp-APAM,Shp-UMB-2018-3-4,Shp-ProcIGC-2018-3-5,Shp-ProcIAMM-2018-3-6}.
This method is systematically described in monograph \cite[chapters 8, 11-13]{Plaksa-Shpakivskyi-2023}.

Let $\mathbb{A}$ be an $n$-dimensional commutative associative algebra
over the field of complex numbers $\mathbb{C}$ and let
$e_1,e_2,e_3\in \mathbb{A}$ be three linearly independent  vectors over the real field $\mathbb{R}$. Let
$\zeta:=xe_1+ye_2+ze_3$, $x,y,z\in\mathbb{R}$.  Let
 $E_3:=\{\zeta=xe_1+ye_2+ze_3:\,\, x,y,z\in\mathbb{R}\}$ be the
linear span of vectors $e_1,e_2,e_3$ over the field
$\mathbb{R}$. Let $\Omega$ be a domain in $E_3$.

\vskip2mm

\textbf{Definition 1.} \textit{
We say that a continuous function
$\Phi:\Omega\rightarrow\mathbb{A}$ is \textit{monogenic}
in $\Omega$ if $\Phi$ is differentiable in the sense of
Gateaux in every point of $\Omega$, i.~e. if  for every
$\zeta\in\Omega$ there exists an element
$\Phi'(\zeta)\in\mathbb{A}$ such that
\begin{equation}\label{monogennaOZNA}\medskip
\lim\limits_{\varepsilon\rightarrow 0+0}
\frac{\Phi(\zeta+\varepsilon
h)-\Phi(\zeta)}{\varepsilon}= h\Phi'(\zeta)\quad\forall\,
h\in E_3.\medskip
\end{equation}
$\Phi'(\zeta)$ is theGateaux derivative of the function
$\Phi$ in the point $\zeta$.}

In the paper \cite{Shp-APAM} it is shown that under certain conditions, which are fulfilled in our case, the monogenic function $\Phi$
 has Gateaux derivatives of all orders $\Phi^r$, $r=1,2,3,...$.
Therefore, if a function $\Phi $  is monogenic in every point of
$\Omega$, then
$$\frac{\partial^{4}\Phi}
{\partial x^4}=
e_1^4\,\Phi^{IV}(\zeta),\quad \frac{\partial^{4}\Phi}
{\partial x^2\partial y^2}=
e_1^2e_2^2\,\Phi^{IV}(\zeta),$$
etc.
Consequently,   due to the equality
$$
\Delta^2 \Phi(\zeta)=\Phi^{IV}(\zeta)(e_1^2+e_2^2+e_3^2)^2=0
$$
every monogenic in
$\Omega $ function $\Phi$
satisfies the equation (\ref{3-bih}) everywhere in
$\Omega $ if and only if
\begin{equation}\label{har-ekb}
(e_1^2+e_2^2+e_3^2)^2=0.
\end{equation}

\vskip2mm

\textbf{Definition 2.} \textit{Equation (\ref{har-ekb}) is called the characteristic equation.}

Instead of the algebra $\mathbb{A}$, we will consider \textit{a sequence} of algebras  $\{\mathbb{E}^n_\rho\}_{n=1}^\infty$. 
For any natural\, $n$, let\, $\mathbb{E}^n_\rho$\, be a commutative associative algebra
over the complex field with a basis\, $\{1,\rho,\rho^2,\ldots,\rho^{n-1}\}$, where\, $\rho^{n}=0$.

Let vectors\, $e_1,e_2,e_3\in\mathbb{E}^n_\rho$\,
have the following decompositions with respect to the basis of
the algebra\, $\mathbb{E}^n_\rho$:
\begin{equation}\label{3-6:e_1_e_2_e_3-k-ro}
 e_1:=\sum\limits_{r=0}^{n-1}k_r\,\rho^r\,,\quad e_2:=\sum\limits_{r=0}^{n-1}m_r\,\rho^r\,,\quad e_3:=\sum\limits_{r=0}^{n-1}g_r\,\rho^r,\quad
 k_r,m_r,g_r\in\mathbb{C}.
\end{equation}
Then $$\zeta=xe_1+ye_2+ze_3=\sum\limits_{r=0}^{n-1}(k_rx+m_ry+g_rz)\,\rho^r=:\sum\limits_{r=0}^{n-1}\xi_r\,\rho^r,$$
where $\xi_0:=k_0x+m_0y+g_0z$ and $\xi_r:=k_rx+m_ry+g_rz$, $r=1,2,...,n-1$.
In this case, the decomposition of resolvent $(t-\zeta)^{-1}$ is of the form
\begin{equation}\label{3-6:rozkl-rezol-A_n^m}
(t-\zeta)^{-1}=\sum\limits_{k=0}^{n-1}A_k\,\rho^k\,\qquad\forall\,\zeta\in E_3 \quad \forall\,t\in\mathbb{C}:\,
t\neq \xi\,,
  \end{equation}
where the coefficients\, $A_k$\, are determined by
the following recurrence relations:
\begin{equation}\label{3-6:-rozkl-rezol-A_n^m}
A_0=\frac{1}{t-\xi_0}\,, \quad A_s=\frac{1}{t-\xi_0}(\xi_sA_0+\xi_{s-1}A_1+\cdots+\xi_1A_{s-1}),
\quad s=1,2,\ldots.
\end{equation}

In particular, we have the following initial functions:
$$A_0=\frac{1}{t-\xi_0}\,, \qquad A_1=\frac{\xi_1}{(t-\xi_0)^2},\qquad A_2=\frac{\xi_2}{(t-\xi_0)^2}+\frac{\xi_1^2}{(t-\xi_0)^3},\quad
$$
$$
A_3=\frac{\xi_3}{(t-\xi_0)^2}+\frac{2\xi_1\xi_2}{(t-\xi_0)^3}+\frac{\xi_1^3}{(t-\xi_0)^4},
$$
$$A_4=\frac{\xi_4}{(t-\xi_0)^2}+\frac{2\xi_1\xi_3+\xi_2^2}{(t-\xi_0)^3}+\frac{3\xi_1^2\xi_2}{(t-\xi_0)^4}+
\frac{\xi_1^4}{(t-\xi)^5}\,,
$$
$$A_5=\frac{\xi_5}{(t-\xi_0)^2}+\frac{2\xi_1\xi_4+2\xi_2\xi_3}{(t-\xi_0)^3}+\frac{3\xi_1^2\xi_3+3\xi_1\xi_2^2}{(t-\xi_0)^4}+
\frac{4\xi_1^3\xi_2}{(t-\xi_0)^5}+\frac{\xi_1^5}{(t-\xi_0)^6}\,,
$$
$$A_6=\frac{\xi_6}{(t-\xi_0)^2}+\frac{\xi_3^2+2\xi_1\xi_5+2\xi_2\xi_4}{(t-\xi_0)^3}+
\frac{\xi_2^3+6\xi_1\xi_2\xi_3+3\xi_1^2\xi_4}{(t-\xi_0)^4}$$
$$+\frac{4\xi_1^3\xi_3+6\xi_1^2\xi_2^2}{(t-\xi_0)^5}+\frac{5\xi_1^4\xi_2}{(t-\xi_0)^6}+\frac{\xi_1^6}{(t-\xi_0)^7}\,,
$$
etc.

Let\, $D$\, be a domain in the complex plane and let\, $\Pi:=\{(x,y,z)\in \mathbb{R}^3 : \xi_0=k_0x+m_0y+g_0z\in D\}$.

As shown in \cite{Shp-ProcIAMM-2018-3-6} (see also subsections 13.2.1, 13.3.2 in \cite{Plaksa-Shpakivskyi-2023}) that 
an infinite-dimensional family of solutions of equation
(\ref{3-bih}) on the set $\Pi$\, is of the following form:
\begin{equation}\label{3-6:sim-1}
 \Biggl\{\frac{1}{2\pi i}\int\limits_\Gamma
F(t)A_k\,dt\Biggr\}_{k=0}^\infty,
 \end{equation}
 where\, $F : D\rightarrow\mathbb{C}$\, may be arbitrary holomorphic functions in the domain\, $D$, $\Gamma$ is a closed Jordan rectifiable curve in $D$ that surround the point
 $\xi_0$, and\, $A_k$\, are defined by the recurrent formulas (\ref{3-6:-rozkl-rezol-A_n^m}), where\, $n$\, can increase to infinity.

Let us write some initial functions, which belong to this
family and are also denoted as\,$\{U_k(x,y,z)\}_{k=0}^\infty$\,.
After substituting the functions\, $A_0, A_1,\dots, A_6$\,
into expressions (\ref{3-6:sim-1}), we get
$$
U_0=F_0(\xi_0)\,,\qquad U_1=\xi_1F'_1(\xi_0)\,, \qquad
U_2=\xi_2F_2'(\xi_0)+\frac{\xi_1^2}{2!}\,F_2''(\xi_0)\,,
$$
$$U_3=\xi_3F_3'(\xi_0)+\xi_1\xi_2F_3''(\xi_0)+\frac{\xi_1^3}{3!}\,F_3'''(\xi_0)\,,\vspace*{2mm}$$
$$U_4=\xi_4F_4'(\xi_0)+\frac{1}{2}(2\xi_1\xi_3+\xi_2^2)F_4''(\xi_0)+\frac{\xi_1^2\xi_2}{2}\,F_4'''(\xi_0)+
\frac{\xi_1^4}{4!}\,F_4^{(4)}(\xi_0)\,,$$
$$
U_5=\xi_5F_5'(\xi_0)+ (\xi_1\xi_4+\xi_2\xi_3)F_5''(\xi_0)+\frac{1}{2}(
\xi_1\xi_2^2+\xi_1^2\xi_3)F_5'''(\xi_0)$$
$$
+\frac{1}{3!}\,\xi_1^3\xi_2F_5^{(4)}(\xi_0)+\frac{\xi_1^5}{5!}\,F_5^{(5)}(\xi_0)\,,
$$
$$
U_6=\xi_6F_6'(\xi_0)+
\frac{1}{2}(\xi_3^2+2\xi_1\xi_5+2\xi_2\xi_4)F_6''(\xi_0)+\frac{1}{6}(
\xi_2^3+6\xi_1\xi_2\xi_3+3\xi_1^2\xi_4)F_6'''(\xi_0)$$
$$
+\frac{1}{4!}( 4\xi_1^3\xi_3+6\xi_1^2\xi_2^2)F_6^{(4)}(\xi_0)+
\frac{\xi_1^4\xi_2}{4!}\,F_6^{(5)}(\xi_0)+\frac{\xi_1^6}{6!}\,F_6^{(6)}(\xi_0)\,,
$$
etc., where\, $F_m$\, for\, $m=0,1,2,3,4,5,6$\, are arbitrary holomorphic
functions of a complex variable.

Calculating the subsequent functions\, $A_k$\, by the recurrent formula
(\ref{3-6:-rozkl-rezol-A_n^m}), we can write the  infinite-dimensional family of solutions
 of equation (\ref{3-bih}) that represent family of form (\ref{3-6:sim-1}).

\section{An infinite-dimensional family of exact solutions of a three-dimensional biharmonic equation}

We shall construct solutions of form (\ref{3-6:sim-1}) for the
three-dimensional biharmonic equation  (\ref{3-bih}). For this purpose, in the sequence of algebras
$\{\mathbb{E}^n_\rho\}_{n=1}^\infty$\,, we shall find all triads of
vectors\, $e_1,e_2,e_3$\, satisfying the characteristic equation (\ref{har-ekb}).

Using decomposition (\ref{3-6:e_1_e_2_e_3-k-ro}) and multiplication table of the algebra   $\mathbb{E}^n_\rho$, we have
\begin{equation}\label{3-6:dec-e^2}
 e_1^2=\sum\limits_{r=0}^{n-1} B_r\rho^r\,,
\end{equation}
where
$$
B_0=k_0^2\,,\quad B_1=2k_0k_1\,,\quad B_2=k_1^2+2k_0k_2\,,
$$
and, in the general case,
$$
B_r\equiv B_r(k_0,k_1,\ldots,k_r):=\left\{
\begin{array}{lrr}
&k_{r/2}^2+2\left(k_0k_r+k_1k_{r-1}+\cdots+k_{\frac{r}{2}-1}k_{\frac{r}{2}+1}\right)\\[1.5mm]
&\;\;\mbox{if}\;\;  r \;\;\mbox{is\, even},\\[1.5mm]
&2\left(k_0k_r+k_1k_{r-1}+\cdots+k_{\frac{r-1}{2}}k_{\frac{r+1}{2}} \right)\\[1.5mm]
&\;\;\mbox{if}\;\;  r \;\;\mbox{is\,odd}.
\end{array}
\right.
$$

In a similar way, we have
\begin{equation}\label{3-6:dec-e2^2}
e_2^2=\sum\limits_{r=0}^{n-1} C_r\rho^r\,,
\end{equation}
where the coefficients\, $C_r$\, satisfies the relations
$$
C_r\equiv C_r(m_0,m_1,\ldots,m_r)=B_r(m_0,m_1,\ldots,m_r)
$$
for\, $r=0,1,\dots,n-1$.

Similarly,
$$
e_3^2=\sum\limits_{r=0}^{n-1} D_r\rho^r\,,
$$
where the coefficients\, $D_r$\, satisfies the relations
$$
D_r\equiv D_r(g_0,g_1,\ldots,g_r)=B_r(g_0,g_1,\ldots,g_r)
$$
for\, $r=0,1,\dots,n-1$.

Now, using (\ref{3-6:dec-e^2}) and multiplication table of the algebra   $\mathbb{E}^n_\rho$, we have
$$
e_1^4=\sum\limits_{r=0}^{n-1} G_r\rho^r\,,
$$
where the coefficients\, $G_r$\, satisfies the relations
\begin{equation}\label{3-6:G_r=B_r}
G_r\equiv G_r(k_0,k_1,\ldots,k_r)=B_r(B_0,B_1,\ldots,B_r)
\end{equation}
for\, $r=0,1,\dots,n-1$.

Similarly,
$$
e_2^4=\sum\limits_{r=0}^{n-1} H_r\rho^r\,,
$$
where the coefficients\, $H_r$\, satisfies the relations
\begin{equation}\label{3-6:H_r=C_r}
H_r\equiv H_r(m_0,m_1,\ldots,m_r)=C_r(C_0,C_1,\ldots,C_r)=G_r(m_0,m_1,\ldots,m_r)
\end{equation}
for\, $r=0,1,\dots,n-1$,
and
$$
e_3^4=\sum\limits_{r=0}^{n-1} M_r\rho^r\,,
$$
where the coefficients\, $M_r$\, satisfies the relations
\begin{equation}\label{3-6:M_r=D_r}
M_r\equiv M_r(g_0,g_1,\ldots,g_r)=D_r(D_0,D_1,\ldots,D_r)=G_r(g_0,g_1,\ldots,g_r)
\end{equation}
for\, $r=0,1,\dots,n-1$.

Equalities (\ref{3-6:dec-e^2}), (\ref{3-6:dec-e2^2}) yield 
$$
e_1^2e_2^2=\sum\limits_{r=0}^{n-1} P_r\rho^r\,,
$$
where the coefficients\, $P_r$\, satisfies the relations
\begin{equation}\label{3-6:P_r=BC_r}
P_r=B_0C_r+B_1C_{r-1}+\cdots+B_rC_0
\end{equation}
for\, $r=0,1,\dots,n-1$.

Similarly,
$$
e_1^2e_3^2=\sum\limits_{r=0}^{n-1} R_r\rho^r\,,
$$
where the coefficients\, $R_r$\, satisfies the relations
\begin{equation}\label{3-6:R_r=BD_r}
R_r=B_0D_r+B_1D_{r-1}+\cdots+B_rD_0
\end{equation}
for\, $r=0,1,\dots,n-1$,
and
$$
e_2^2e_3^2=\sum\limits_{r=0}^{n-1} S_r\rho^r\,,
$$
where the coefficients\, $S_r$\, satisfies the relations
\begin{equation}\label{3-6:S_r=CD_r}
S_r=C_0D_r+C_1D_{r-1}+\cdots+C_rD_0
\end{equation}
for\, $r=0,1,\dots,n-1$.

Thus, characteristic equation (\ref{har-ekb}) is equivalent to the following system of equations
\begin{equation}\label{har-syst}
G_r+H_r+M_r+2P_r+2R_r+2S_r=0,\qquad r=0,1,2,\ldots ,
\end{equation}
where $G_r,H_r,M_r,P_r,R_r,S_r$ are defined by equalities (\ref{3-6:G_r=B_r}), (\ref{3-6:H_r=C_r}),
(\ref{3-6:M_r=D_r}), (\ref{3-6:P_r=BC_r}), (\ref{3-6:R_r=BD_r}),  (\ref{3-6:S_r=CD_r}), respectively.

According to Remark 13.1 \cite{Plaksa-Shpakivskyi-2023} in characteristic equation (\ref{har-ekb}) we can set the vectors\, $e_1, e_2$\, arbitrarily, i.e.,\,
$k_r$\, and\, $m_r$\,  may be arbitrary complex numbers for\, $r=0,1,\ldots, n-1$\, for any natural\, $n$.
Then the vector\, $e_3$\, can be found, and the coefficients\, $g_r$\, can be found from system  (\ref{har-syst}).

We note that for every fixed $r$ equation (\ref{har-syst}) is linear with respect to $g_r$. Therefore, it is easy to find $g_r$  for each fixed $r=0,1,2,\ldots$.

So, characteristic equation (\ref{har-ekb}) is solved on the sequence of algebras $\{\mathbb{E}^n_\rho\}_{n=1}^\infty$.

Thus, we proved the following theorem

\vskip2mm

\textbf{Theorem.} \textit{In a domain $\Pi$ equation (\ref{3-bih}) has 
infinite-dimensional family of exact solutions  of the  form
(\ref{3-6:sim-1}),
  where\, $F : D\rightarrow\mathbb{C}$\, may be arbitrary holomorphic functions in the domain\, $D$, $\Gamma$ is a closed Jordan rectifiable curve in $D$ that surround the point
 $\xi_0$, and\, $A_k$\, are defined by the recurrent formulas (\ref{3-6:-rozkl-rezol-A_n^m}), where\, $n$\, can increase to infinity.  Moreover,
in coefficients $A_k$ the variables $\xi_0$, $\xi_r$ are of the form $\xi_0 =k_0x+m_0y+g_0z$, 
 $\xi_r=k_rx+m_ry+g_rz$, $r=1,2,\ldots$, where $k_r,m_r$ for $r=0,1,2,\ldots$, are arbitrary
complex numbers, and $g_r$ for each fixed $r=0,1,2,\ldots$ are determined by equality  (\ref{har-syst}).
}

\section{Several first solutions from family   (\ref{3-6:sim-1})}

Now we will write some first solutions from   family (\ref{3-6:sim-1}).
For $r=0$ equation (\ref{har-syst}) has the form
\begin{equation}\label{g_0}
(k_0^2+m_0^2+g_0^2)^2=0 \qquad {\rm or} \qquad 
g_0=\pm i\sqrt{k_0^2+m_0^2}\,.
\end{equation}
Therefore, $\xi_0=k_0x+m_0y\pm i\sqrt{k_0^2+m_0^2}\,z$, where $k_0,m_0$ are arbitrary
complex numbers.

Thus, the \textbf{first exact solution} of equation (\ref{3-bih})  from   family (\ref{3-6:sim-1}) is of the form
$$
U_0(x,y,z)=F_0(\xi_0)=F_0\left(k_0x+m_0y\pm i\sqrt{k_0^2+m_0^2}\,z\right),
$$
where $F_0$ is an arbitrary holomorphic
function of the complex variable  $\xi_0=k_0x+m_0y\pm i\sqrt{k_0^2+m_0^2}\,z$ with arbitrary
complex numbers $k_0,m_0$.

For $r=1$ equation (\ref{har-syst}) has the following form
$$
2k_0^3k_1+2m_0^3m_1+2g_0^3g_1+4k_0^2m_0m_1+4m_0^2k_0k_1
$$
\begin{equation}\label{g_1-0}
+4k_0^2g_0g_1+4k_0k_1g_0^2+4m_0^2g_0g_1+4m_0m_1g_0^2=0,
\end{equation}
i.e., we have the linear equation with respect to $g_1$.
Substituting expression (\ref{g_0}) into (\ref{g_1-0}), we obtain
\begin{equation}\label{g_1}
g_1=\pm i\frac{k_0^3k_1+m_0^3m_1}{2\left(\sqrt{k_0^2+m_0^2}\right)^3}\,.
\end{equation}
Therefore, 
$$
\xi_1=k_1x+m_1y \pm i\frac{k_0^3k_1+m_0^3m_1}{2\left(\sqrt{k_0^2+m_0^2}\right)^3}\,z\,,
$$
where $k_0,m_0,k_1,m_1$ are arbitrary
complex numbers.

Thus, the \textbf{second exact solution} of equation (\ref{3-bih})  from   family (\ref{3-6:sim-1}) is of the form
$$
U_1(x,y,z)=\xi_1F_1(\xi_0)$$
$$
=\left(k_1x+m_1y \pm i\frac{k_0^3k_1+m_0^3m_1}{2\left(\sqrt{k_0^2+m_0^2}\right)^3}\,z\right)
F_1\left(k_0x+m_0y\pm i\sqrt{k_0^2+m_0^2}\,z\right),
$$
where $F_1$ is an arbitrary holomorphic
function of the complex variable  $\xi_0=k_0x+m_0y\pm i\sqrt{k_0^2+m_0^2}\,z$, and $k_0,m_0,k_1,m_1$ are arbitrary
complex numbers.

For $r=2$ equation (\ref{har-syst}) has the following form
$$
4k_0^6k_1^2+4k_0^4k_2+2k_0^2k_1^2+4m_0^6m_1^2+4m_0^4m_2+2m_0^2m_1^2+4g_0^6g_1^2+4g_0^4g_2+2g_0^2g_1^2
$$
$$
+2k_0^2m_1^2+4k_0^2m_0m_2+8k_0k_1m_0m_1+2m_0^2k_1^2+4m_0^2k_0k_2+2k_0^2g_1^2+4k_0^2g_0g_2
$$
$$+8k_0k_1g_0g_1+2k_1^2g_0^2+4k_0k_2g_0^2+2m_0^2g_1^2+4m_0^2g_0g_2+8m_0m_1g_0g_1
$$
\begin{equation}\label{g_2-0}
+2m_1^2g_0^2+4m_0m_2g_0^2=0.
\end{equation}
i.e., we have the linear equation with respect to $g_2$.
From (\ref{g_2-0}) we obtain
$$
g_2=\frac{1}{4(g_0^3-g_0^4)}\Biggl(
4k_0^6k_1^2+4k_0^4k_2+2k_0^2k_1^2+4m_0^6m_1^2+4m_0^4m_2+2m_0^2m_1^2
$$
$$
+4g_0^6g_1^2+2g_0^2g_1^2+2k_0^2m_1^2+4k_0^2m_0m_2+8k_0k_1m_0m_1+2m_0^2k_1^2+4m_0^2k_0k_2
$$
$$+2k_0^2g_1^2+8k_0k_1g_0g_1+2k_1^2g_0^2+4k_0k_2g_0^2+2m_0^2g_1^2+8m_0m_1g_0g_1
$$
\begin{equation}\label{g_2}
+2m_1^2g_0^2+4m_0m_2g_0^2\Biggr),
\end{equation}
where $g_0,g_1$ are defined by (\ref{g_0}), (\ref{g_1}), and $k_0,m_0,k_1,m_1,k_2,m_2$ are arbitrary
complex numbers.

Therefore, 
$$
\xi_2=k_2x+m_2y+ \frac{1}{4(g_0^3-g_0^4)}\Biggl(
4k_0^6k_1^2+4k_0^4k_2+2k_0^2k_1^2+4m_0^6m_1^2+4m_0^4m_2
$$
$$
+2m_0^2m_1^2+4g_0^6g_1^2+2g_0^2g_1^2+2k_0^2m_1^2+4k_0^2m_0m_2+8k_0k_1m_0m_1+2m_0^2k_1^2
$$
$$+4m_0^2k_0k_2+2k_0^2g_1^2+8k_0k_1g_0g_1+2k_1^2g_0^2+4k_0k_2g_0^2+2m_0^2g_1^2
$$
$$
+8m_0m_1g_0g_1+2m_1^2g_0^2+4m_0m_2g_0^2\Biggr)\,z\,,
$$
where $k_0,m_0,k_1,m_1,k_2,m_2$ are arbitrary
complex numbers, and $g_0,g_1$ are defined by equalities (\ref{g_0}), (\ref{g_1}).

Thus, the \textbf{third exact solution} of equation (\ref{3-bih})  from   family (\ref{3-6:sim-1}) is of the form
$$
U_2=\xi_2F_2(\xi_0)+\frac{\xi_1^2}{2!}\,F_2'(\xi_0)$$
$$
=\left[k_2x+m_2y+ \frac{1}{4(g_0^3-g_0^4)}\Biggl(
4k_0^6k_1^2+4k_0^4k_2+2k_0^2k_1^2+4m_0^6m_1^2+4m_0^4m_2\right.
$$
$$
+2m_0^2m_1^2+4g_0^6g_1^2+2g_0^2g_1^2+2k_0^2m_1^2+4k_0^2m_0m_2+8k_0k_1m_0m_1+2m_0^2k_1^2
$$
$$+4m_0^2k_0k_2+2k_0^2g_1^2+8k_0k_1g_0g_1+2k_1^2g_0^2+4k_0k_2g_0^2+2m_0^2g_1^2
$$
$$\left.
+8m_0m_1g_0g_1+2m_1^2g_0^2+4m_0m_2g_0^2\Biggr)\,z\,\right]\,F_2\left(k_0x+m_0y\pm i\sqrt{k_0^2+m_0^2}\,z\right)
$$
$$
+\frac{1}{2}\left(k_1x+m_1y \pm i\frac{k_0^3k_1+m_0^3m_1}{2\left(\sqrt{k_0^2+m_0^2}\right)^3}\,z\right)^2
F_2'\left(k_0x+m_0y\pm i\sqrt{k_0^2+m_0^2}\,z\right),
$$
where $F_2$ is an arbitrary holomorphic
function of the complex variable  $\xi_0=k_0x+m_0y\pm i\sqrt{k_0^2+m_0^2}\,z$, and $k_0,m_0,k_1,m_1,k_2,m_2$ are arbitrary
complex numbers, and $g_0,g_1$ are defined by equalities (\ref{g_0}), (\ref{g_1}).

Continuing in this way we can write out an infinite number of exact solutions of equation (\ref{3-bih}).

\vskip2mm
\textbf{The algorithm.}
\begin{enumerate}
  \item 
  
  \begin{itemize}
          \item We put $r=0$;
            \item from the linear equation (\ref{har-syst}) we find $g_0$ and we write out $\xi_0=k_0x+m_0y+g_0z$. Here $k_0,m_0$ are arbitrary
complex numbers;            
  \item $\xi_0$ we substitute in $A_0$ from (\ref{3-6:-rozkl-rezol-A_n^m});
  \item $A_0$ we substitute in  (\ref{3-6:sim-1}) and we obtain the solution $U_0$.
        \end{itemize}
        
    \item 
    
    \begin{itemize}
          \item We put $r=1$;
            \item having $g_0$ from the linear equation (\ref{har-syst}) we find $g_1$ and we write out $\xi_1=k_1x+m_1y+g_1z$. 
            Here $k_0,m_0,k_1,m_1$ are arbitrary
complex numbers;            
  \item $\xi_0$ and $\xi_1$ we substitute in $A_1$ from (\ref{3-6:-rozkl-rezol-A_n^m});
  \item $A_1$ we substitute in  (\ref{3-6:sim-1}) and we obtain the solution $U_1$.
        \end{itemize}
        
        \item 
    
    \begin{itemize}
          \item We put $r=2$;
            \item having $g_0,g_1$ from the linear equation (\ref{har-syst}) we find $g_2$ and we write out $\xi_2=k_1x+m_1y+g_1z$.
             Here $k_0,m_0,k_1,m_1,k_2,m_2$ are arbitrary
complex numbers;
  \item $\xi_0$, $\xi_1$ and $\xi_2$ we substitute in $A_2$ from (\ref{3-6:-rozkl-rezol-A_n^m});
  \item $A_2$ we substitute in  (\ref{3-6:sim-1}) and we obtain the solution $U_2$.
        \end{itemize}
        
       \item   etc.
        
\end{enumerate}

 \subsection*{Acknowledgment}
This work was supported by a grant from the Simons Foundation
(1030291,V.S.Sh.).

\renewcommand{\refname}{References}

\vskip10mm
\noindent Vitalii Shpakivskyi

\noindent Institute of Mathematics of the\\
National Academy of Sciences of Ukraine, Kyiv

\noindent e-mail: shpakivskyi86@gmail.com

\end{document}